# Davenport constant for finite abelian groups with higher rank

by Anamitro Biswas [A], Eshita Mazumdar[B]

**Abstract**

For a finite abelian group $G$, the Davenport Constant, denoted by $D(G)$, is defined to be the least positive integer $k$ such that every sequence of length at least $k$ has a non-trivial zero-sum subsequence. A long-standing conjecture is that the Davenport constant of a finite abelian group $G = C_{n_1} \times \cdots \times C_{n_d}$ of rank $d \in \mathbb{N}$ is $1 + \sum_{i=1}^{d} (n_i - 1)$. This conjecture is false in general, but it remains to know for which groups it is true. In this paper, we consider groups of the form $G = (C_p)^{d-1} \times C_{pq}$, where $p$ is a prime and $q \in \mathbb{N}$ and provide sufficient condition when the conjecture holds true.

## 1 Introduction

Our notations are mostly standard. For integers $a, b$ we set $[a, b] = \{x \in \mathbb{Z} : a \le x \le b\}$. Let $G$ be a finite abelian group written additively. By a $G$-sequence, we shall mean a finite sequence $x_1 \ldots x_k$ with $x_i \in G$ for each $i$, where repetition of element is allowed. If we write, a sequence $S = x_1^{n_1} \cdots x_s^{n_s}$ with $x_i \ne x_j$ for $i \ne j$, it means in the sequence $S$ the element $x_i$ appears $n_i$ times, and $n_i = v_{x_i}(S)$ is said to be the multiplicity of $x_i$ for each $i$. A sequence $T$ is said to be a subsequence of $S$, denoted by $T \,|\, S$, if $v_g(T) \le v_g(S)$ for every $g$ in $S$ and $ST^{-1}$ denotes the subsequence obtained from $S$ by removing elements of $T$. Two subsequences $T$ and $U$ of $S$ are said to be disjoint subsequences of $S$ if $U$ is a subsequence of $S\,T^{-1}$.

By a zero-sum $G$-sequence, we shall mean a $G$-sequence $x_1 \ldots x_k$ such that $\sum_i x_i = 0_G$, where $0_G$ is the identity element of $G$. The Davenport Constant $D(G)$, introduced by Rogers [R], is defined as the smallest integer $k$ such that every $G$-sequence of length $k$ contains a non-trivial zero-sum subsequence. Its precise value for any finite abelian group is yet to be known. The Davenport constant is an important invariant of the ideal class group of the ring of integers of an algebraic number field (see [H] for more details).

We know that every finite abelian group $G$ can be written as a product of cyclic groups as $C_{n_1} \times \cdots \times C_{n_d}$ where $n_i \in \mathbb{N}$ with $n_i \,|\, n_{i+1}$. Then the integer $n_d$ is said to be the *exponent* of the group $G$, denoted by $\exp(G)$, and $d$ is said to be the *rank* of the group $G$. For $d = 1, 2$ we have $D(G) = n_1$, and $D(G) = n_1 + n_2 - 1$ ([O2]) respectively. Let us consider $D^*(G) = 1 + \sum_{i=1}^{d}(n_i - 1)$. Clearly, $D(G) \ge D^*(G)$. A fairly long-standing conjecture by Olson is that $D(G) = D^*(G)$. He proved the conjecture for finite abelian $p$-groups for prime $p$ in [O]. This conjecture is false in general, but it remains to know for which groups this is true. In [GeS] Geroldinger and Schneider proved that there are infinitely many groups $G$ of rank $4, 5$ with $D(G) > D^*(G)$. Also, the problem of existence of a group $G$ of rank 3 with $D(G) > D^*(G)$ remains open.

*A.* Department of Mathematics, NIT Silchar, Assam -788010, India. email: anamitroappu(at)gmail.com

*B.* Mathematical and Physical Sciences, School of Arts and Sciences, Ahmedabad University, Gujrat - 380009, India. email: eshita.mazumdar(at)ahduni.edu.in As of now, it is hard to say anything specifically about groups of rank 3 or more. In this paper, we consider a certain class of groups of rank $d \geq 3$ i.e., the group of the form $G = C_p^{d-1} \times C_{pq}$, for prime $p$, and $q \in \mathbb{N}$ and we shall investigate when $D(G) = D^*(G)$. Such type of work has been done in [S] for the group of the type $C_p \times C_p \times C_{2p}$, where $p \neq 2$ is a prime. Now, we have the result in more general setting. We tackle such problems in two ways, one way includes $(p, q) = 1$ and other way includes $(p, q) \neq 1$. We shall prove that $D(G) = D^*(G)$, for $(p, q) = 1$ under certain sufficient condition and we also show that for $(p, q) \neq 1$, $D(G) = D^*(G)$ follows directly. First, we shall be stating the sufficient condition in the form of a following conjecture:

**Conjecture 1.** *Let us consider prime $p$, and $q \in \mathbb{N}$ such that $(p, q) = 1$. Let us define $G := C_p^d \times C_q$, where $d \in \mathbb{N}^{\geq 3}$. Also, assume that $x_1 \dots x_m$ is a sequence over $C_p^d$ and $y_1 \dots y_m$ is a sequence over $C_q$ where $m = p(q + d - 1) - (d - 1)$. Suppose*

$$y_1 = \dots = y_{r_1} = 1,$$
$$y_{r_1+1} = \dots = y_{r_1+r_2} = 2,$$
$$\vdots$$
$$y_{\sum_{i=1}^{q-2} r_i + 1} = \dots = y_{\sum_{i=1}^{q-1} r_i} = q - 1,$$
$$y_{r+1} = \dots = y_m = 0 \text{ where } r = \sum_{i=1}^{q-1} r_i.$$

*If $r \in [pq + 1, p(q + d - 1) - (d - 1)]$ and $\sum_{i=1}^{q-1} i \, r_i \equiv 0 \, (mod \, q)$, then the sequence $S := (x_1, y_1) \dots (x_m, y_m)$ over $G$ has a nontrivial zero-sum subsequence.*

It is already known that $D(G) = D^*(G)$ for the groups $C_3 \times C_3 \times C_{3s}$ for any positive integer $s$ ([BP2]). So, Conjecture 1 is true for $p = 3$ with $d = 3$ at least. Also, from the above conjecture we can conclude our main theorem.

**Theorem 2.** *Let $p$ be a prime and $q \in \mathbb{N}$ with $(p, q) = 1$ such that Conjecture 1 holds. Then, for group $G = C_p^{d-1} \times C_{pq}$ with $d \in \mathbb{N}^{\geq 3}$, $D(G) = D^*(G)$.*

For dealing the second part i.e, $(p, q) \neq 1$, we need to discuss some properties of a generalized version of Davenport constant which is known as $r$-wise Davenport constant ([FS]), which was introduced by Halter-Koch in 1992 ([HK]). For a finite abelian group $G$ and $r \in \mathbb{N}$, $r$-wise Davenport constant, denoted by $D_r(G)$, is defined as the least positive integer $k$ such that every sequence of length at least $k$ has $r$ disjoint zero-sum subsequences. Clearly, $D_r(G) \leq D_{r+1}(G)$. Various other information about this constant will be available in [BP], [DOQ], [GZZ]. For finite abelian group of rank $d = 1, 2$, we have $D_r(G) = r \, n_1$, and $D_r(G) = r \, n_2 + n_1 - 1$ respectively ([GS]). Now it's natural to ask the precise value of $D_r(G)$ for any finite abelian group of rank $d \geq 3$. People dealt similar problem for certain finite abelian $p$-group for prime $p$ (Theorem 6.1.5, [GH]). One can determine $D_r(G)$ for that particular class. In this connection, we have the value of $r - $ wise Davenport constant for finite abelian $p$-group for prime $p$. There is an approach mentioned in Theorem 6.1.5 of [GH], from which one can determine the same for certain class of $p$-groups. We also followed slightly different approach to prove the same result. The result for $p$-groups is the following:

**Theorem 3.** *Let $p$ be a prime and $e_i, d \in \mathbb{N}$ for $i \in [1, d]$. For a finite abelian $p$-group $G = C_{p^{e_1}} \times C_{p^{e_2}} \times \dots \times C_{p^{e_d}}$ such that $e_i \leq e_{i+1}$, $1 \leq i \leq d - 1$ and $p^{e_d} \geq 1 + \sum_{i=1}^{d-1} (p^{e_i} - 1)$, we have $D_r(G) = r \, p^{e_d} + \sum_{i=1}^{d-1} p^{e_i} - d + 1$.*

**Theorem 4.** *Let $p$ be an odd prime and $e_i \in \mathbb{N}$ for $i \in \{1, 2, \ldots, d\}$. Then for the group $G = C_{p^{e_1}} \times C_{p^{e_2}} \times \ldots \times C_{p^{e_{d-1}}} \times C_{mp^{e_d}}$ with $e_1 \leq e_2 \leq e_3 \leq \ldots \leq e_d$ such that $p^{e_d} \geq 1 + \sum_{i=1}^{d-1} (p^{e_i} - 1)$, we have,*

$$D_r(G) = r \, n \, p^{e_d} + \sum_{i=1}^{d-1} p^{e_i} - d + 1.$$

**Corollary 5.** *Let $p$ be an odd prime and $G = C_{p^{e_1}} \times C_{p^{e_2}} \times \ldots \times C_{p^{e_{d-2}}} \times C_{mp^{e_{d-1}}} \times C_{np^{e_d}}$ be a group with $e_1 \leq e_2 \leq e_3 \leq \ldots \leq e_d$ and $p^{e_d} \geq 1 + \sum_{i=1}^{d-1} (p^{e_i} - 1)$.*

1. *If $m \,|\, n$ then,*

$$r \, n \, p^{e_d} + (m-1) \, p^{e_{d-1}} + \sum_{i=1}^{d-1} p^{e_i} - d + 1 \leq D_r(G) \leq r \, n \, p^{e_d} + (m-1) \, p^{e_d} + \sum_{i=1}^{d-1} p^{e_i} - d + 1.$$

2. *If $n \,|\, m$ then,*

$$\max \left\{ r \, n \, p^{e_d} + (m-1) \, p^{e_{d-1}} + \sum_{i=1}^{d-1} p^{e_i} - d + 1, \, n \, p^{e_d} + (r \, m - 1) \, p^{e_{d-1}} + \sum_{i=1}^{d-1} p^{e_i} - d + 1 \right\}$$

$$\leq \quad D_r(G)$$

$$\leq \quad r \, m \, p^{e_d} + (n-1) \, p^{e_d} + \sum_{i=1}^{d-1} p^{e_i} - d + 1.$$

From the above result, we conclude the following bounds:

**Theorem 6.** *Consider the group $G = C_{\prod_{j=1}^{\ell} p_1^{e_1^{(j)}}} \times \ldots \times C_{\prod_{j=1}^{\ell} p_d^{e_d^{(j)}}}$, where w.l.o.g we can assume that $e_1^{(j)}, e_2^{(j)}, \ldots, e_d^{(j)}$ $(1 \leq j \leq \ell)$ are integers satisfying $0 \leq e_1^{(j)} \leq e_2^{(j)} \leq \ldots \leq e_d^{(j)}$ for all $j$ but for each $j \in [1, \ell]$ all $e_i^{(j)}$'s are not zero for $i \in [1, d]$. If $p_j$'s are distinct primes such that $p_j^{e_d^{(j)}} \geq 1 + \sum_{i=1}^{d-1} \left( p_j^{e_i^{(j)}} - 1 \right)$ for all $j = [1, \ell]$. Then, for $\varphi(p_j) = \sum_{i=1}^{d-1} p_j^{e_i^{(j)}} - d + 1$ with $j = [1, \ell]$, we have*

$$r \prod_{j=1}^{\ell} p_j^{e_d^{(j)}} + \sum_{i=1}^{d-1} \left( \prod_{j=1}^{\ell} p_j^{e_i^{(j)}} - 1 \right) \quad \leq \quad D_r(G)$$

$$\leq \quad r \prod_{j=1}^{\ell} p_j^{e_d^{(j)}} + \sum_{m=1}^{\ell-1} \left( \left( \prod_{j=m+1}^{\ell} p_j^{e_d^{(j)}} \right) \varphi(p_m) \right) + \varphi(p_\ell),$$

**Observation 1:** It follows from theorem 6 that if we have the class of groups of the form $G = C_{p^{e_1}} \times C_{p^{e_2}} \times \ldots \times C_{p^{e_{d-1}}} \times C_{mp^{e_d}}$ with $e_i \in \mathbb{N}$, and $e_i \leq e_{i+1}$ such that $p^{e_d} \geq 1 + \sum_{i=1}^{d-1} (p^{e_i} - 1)$, then $D_r(G) = r \, m \, p^{e_d} + \sum_{i=1}^{d-1} p^{e_i} - d + 1$.

We conclude the second part of our main theorem from Observation 1:

**Theorem 7.** *For a sufficiently large prime $p$ and $q \in \mathbb{N}$ with $(p, q) \neq 1$. Then, for group $G = C_p^{d-1} \times C_{pq}$ with $d \in \mathbb{N}^{\geq 3}$, $D(G) = D^*(G)$.*

**Observation 2:** From theorem 6, one can observe that the magnitude of the ratio of the upper and lower bounds of $D_r(G)$ which will tend to one for either higher powers of primes, or sufficiently large primes, or sufficiently large values of $r$. We shall examine it in the section 6.

Most interestingly it could be noted that though we are able to conclude about $D_r(G)$ for certain class of $G$ with $\exp(G) = p^{e_d}$, in theorem 3 but we are not able to conclude anything about it if $\exp(G) = p$ with $d \geq 3$. We have observed that $D_r((C_p)^d) = (r + d - 1) p - (d - 1)$ holds for groups of rank upto 2 ([GS]). Also the same is true for $r = 1$ ([O]). But it is hard to find a general form of $D_r((C_p)^d)$. The rest of the paper is organised as follows: we begin with some preliminaries before launching into some useful lemmas in section 3. In section 4, we prove theorem 2. In next section we will be dealing with our theorem for the case $(p, q) \neq 1$ and prove theorems related to $r$-wise Davenport constant. We conclude the paper with remarks.

⟨iffalse⟩ Then in section 3, we will discuss how Conjecture 1 will imply Theorem 2. In section 4, we will prove how Conjecture ? will imply Theorem ? and how Conjecture ? will imply Conjecture 1. Finally we will discuss some concluding remarks in section 5. ⟨fi⟩

## 2   Preliminaries

If $G$ is a finite abelian group and $S$ is a sequence over $G$ then

- $|S|$, $\sigma(S)$, and $[S]$ denotes the length of the sequence, sum of all the elements of $S$, and the set of all sums of non-trivial subsequences of $S$, respectively.

- $S$ is said to zero-sum free, and a short zero-sum sequence if it has no nonempty zero-sum $G$-subsequence, and if it is a zero-sum $G$-sequence of length $|S| \in [1, e\,x\,p(G)]$, respectively.

- If $G = H_1 \times H_2$ where $H_i$'s are finite abelian groups for $i = 1, 2$. Let $S = (x_1, y_1) \ldots (x_k, y_k)$ be a sequence over $G$. Given any non-empty subsequence $T = x_{i_1} \ldots x_{i_s}$ of $x_1 x_2 \ldots x_k$, we define the extension of $T$ into $S$, denoted by $T^S$ to be the subsequence $(x_{i_1}, y_{i_1}) \ldots (x_{i_s}, y_{i_s})$.

We can recall the definition of the zero-sum invariant $\eta(G)$ to be the smallest positive integer $k$ such that every $G$-sequence of length at least $k$ has a short zero-sum subsequence. Clearly $D(G) \leq \eta(G)$. In this paper we introduce a new invariant $\eta_r(G)$, which we define as the smallest positive integer $k$ such that every $G$-sequence of length at least $k$ has $r$ disjoint short zero-sum subsequences. Clearly, for $r = 1$, $\eta_1(G) = \eta(G)$. So $\eta_r(G)$ is a natural generalization of $\eta(G)$. We will discuss some interesting results about $\eta_r(G)$ in subsequent sections. In this section, we state a few well known results that we serve as important tools in the proof our results: It has been observed by Delorme et. al. [DOQ] that there exists a relationship between $D_r(G)$ and $D_m(G/H)$ where $H$ is a normal subgroup of $G$. We are going to use the following relationship throughout the paper for establishing our results: ⟨fi⟩

**Lemma 8.** *[[DOQ], Proposition 2.6] Let $H$ be a normal subgroup of a finite group $G$ and $r \in \mathbb{N}$, then $D_r(G) \leq D_{D_r(H)}(G/H)$ and $D(G) \geq D(H) + D(G/H) - 1$.*

**Lemma 9.** *[[FGWZ], Lemma 2.7] Let $m$ be a positive integer, and $H$ be a finite abelian group with $\exp(H)|m$ and $m \geq D(H)$. Suppose that $D(C_m \times C_m \times H) = 2\,m + D(H) - 2$. Then $\eta_1(C_m \times H) \leq 2\,m + D(H) - 2$.*

## 3 Some useful results

**Lemma 10.** *Let $p$ be a prime number. Consider a p-group $G = C_{p^{e_1}} \times C_{p^{e_2}} \times \ldots \times C_{p^{e_d}}$ with $1 \leq e_i \leq e_{i+1}$ for $i \in [1, d-1]$. If $p^{e_d} \geq 1 + \sum_{i=1}^{d-1} (p^{e_i} - 1)$, then $\eta(G) \leq D(G) + exp(G)$.*

**Proof.** For the given group $G = C_{p^{e_1}} \times C_{p^{e_2}} \times \ldots \times C_{p^{e_d}}$, let us consider $H = C_{p^{e_1}} \times C_{p^{e_2}} \times \ldots \times C_{p^{e_{d-1}}}$. Clearly, $\exp(H) = p^{e_{d-1}}$ and $D(H) = 1 + \sum_{i=1}^{d-1} (p^{e_i} - 1)$. Let us assume that $m = p^{e_d}$. Now we can see that $\exp(H) \mid m$ and $m \geq D(H)$ (from the given condition). We have, $D\left(C_{p^{e_d}} \times C_{p^{e_d}} \times H\right) = 2\, p^{e_d} + p^{e_{d-1}} + \ldots + p^{e_1} - d$. On the other hand, $2\, m + D(H) - 2 = 2\, p^{e_d} + p^{e_{d-1}} + \ldots + p^{e_1} - d$. Since $D\left(C_{p^{e_d}} \times C_{p^{e_d}} \times H\right) = 2\, p^{e_d} + D(H) - 2$, we can conclude the proof from Lemma 9. $\qquad\square$

**Lemma 11.** *Let $G$ be any finite abelian group and $r \geq 1$ be an integer. If $\eta(G) \leq D(G) + exp(G)$, then $\eta_r(G) \leq D(G) + r\, exp(G)$.*

**Proof.** The lemma is true for $r = 1$. Assume the lemma is true for $r - 1$, with $r \geq 2$. Therefore, $\eta_{r-1}(G) \leq D(G) + (r-1)\exp(G)$. Consider a sequence $S$ over $G$ with $|S| = D(G) + r \exp(G)$. So, $S$ has $(r-1)$ disjoint short zero-sum sub-sequences $S_i$ for $i \in [1, r-1]$. So, $|S\, S_1^{-1}\, S_2^{-1} \ldots S_{r-1}^{-1}| \geq \eta(G)$. $\qquad\square$

**Lemma 12.** *For any finite abelian group $G$, if $\eta_{r-1}(G) \leq D(G) + (r-1) \exp(G)$ then $D_r(G) \leq D(G) + (r-1) \exp(G)$.*

**Proof.** Let $S$ be a sequence over $G$ with $|S| = D(G) + (r-1) \exp(G)$. Since $|S| \geq \eta_{r-1}(G)$ so $S$ has $(r-1)$ disjoint short subsequences $S_i$ for $i \in [1, r-1]$. Therefore, $|S\, S_1^{-1}\, S_2^{-1} \ldots S_{r-1}^{-1}| \geq D(G)$. This completes the proof. $\qquad\square$

**Lemma 13.** *Let $H$ be a finite abelian group and $G = H \times C_q$, where $q \in \mathbb{N}$. Let $S = (x_1, y_1) \ldots (x_\ell, y_\ell)$ be a sequence over $G$. Then $0_G \in [S]$ if either*

*(i). the sequence $x_1 \ldots x_\ell$ contains $q$ non-empty disjoint zero-sum subsequences, or*

*(ii). there exists a non-empty zero-sum subsequence $T$ of $x_1 \ldots x_\ell$ such that $q$ divides $|T|$ and $y_1 = \ldots = y_\ell$,*

**Proof.** Case-I: Assume that the sequence $U = x_1 \ldots x_\ell$ has $q$ non-empty disjoint subsequences $T_i = x_{j_1^{(i)}}\, x_{j_2^{(i)}} \ldots x_{j_{v_i}^{(i)}}$ for $1 \leq i \leq q$ such that $\sigma(T_i) = 0 \; \forall\, i$. Let us consider $q$ non-empty disjoint subsequences $T_i' = y_{j_1^{(i)}}\, y_{j_2^{(i)}} \ldots y_{j_{v_i}^{(i)}}$ for $1 \leq i \leq q$ of $T' = y_1 \ldots y_\ell$ such that $T_i^S = (x_{j_1^{(i)}}, y_{j_1^{(i)}}) \cdots (x_{j_{v_i}^{(i)}}, y_{j_{v_i}^{(i)}})$. Assume that $\sigma(T_i') = y_i' \in C_q$. Since $D(C_q) = q$, so there exists $y_{i_j}'$ for $1 \leq i_1 < i_2 < \ldots < i_m \leq q$ such that $\sigma\left(y_{i_1}'\, y_{i_2}' \ldots y_{i_m}'\right) = 0_{C_q}$, which implies $\sigma\left(T_{i_1}'\, T_{i_2}' \ldots T_{i_m}'\right) = 0_{C_q}$. Therefore, $\sigma\!\left((T_{i_1}\, T_{i_2} \ldots T_{i_m})^S\right) = 0_G$.

Case-II: Consider $\sigma(T) = 0_H$. Let $T = x_{i_1} \ldots x_{i_t}$, where $q \mid t$. Then, $\sigma(T^S) = (0_H, t\, y_1) = (0_H, 0_{C_q}) = 0_G$ as $q \mid t$. $\qquad\square$

**Proposition 14.** *For prime $p$ and $q \in \mathbb{N}$ with $(p, q) = 1$, consider the group $G = C_p^d \times C_q$ with $d \in \mathbb{N}$. Let $x_1 \cdots x_m$ be a sequence over $C_p^d$ and $y_1 \cdots y_m$ be a sequence over $C_q$ where $m = p\,(q + d - 1) - (d - 1)$. W.l.o.g we may assume that for some $r \in [0, m]$,*

$$y_{r+1} = \ldots = y_m = 0.$$

*Then for the sequence $S = (x_1, y_1) \ldots (x_m, y_m)$ over $G$, $r \in [0, p\,q]$ implies $0_G \in [S]$.*

**Proof.** We can assume w.l.o.g that there are $r_i \in \mathbb{N} \cup \{0\}$ for $i \in [1, q-1]$ such that

$$y_1 = \ldots = y_{r_1} = 1,$$
$$y_{r_1+1} = \ldots = y_{r_1+r_2} = 2,$$
$$\vdots$$
$$y_{\sum_{i=1}^{q-2} r_i + 1} = \ldots = y_{\sum_{i=1}^{q-1} r_i} = q-1,$$
$$y_{r+1} = \ldots = y_m = 0 \text{ where } r = \sum_{i=1}^{q-1} r_i.$$

Now consider the sequence $T = \prod_{i=1}^{(q-1)} \prod_{j=1}^{t_i} \sigma(X_i^j)\, x_{r+1}\, x_{r+2} \cdots x_m$, where $t_i = \lfloor \frac{r_i}{q} \rfloor$ for $i \in [1, q-1]$ with $r_i = q\,t_i + \ell_i$, where $l_i \in [0, q-1]$. Note that $X_i^j = x_{\sum_{k=1}^{i-1} r_k + (j-1)q+1} \cdots x_{\sum_{k=1}^{i-1} r_k + jq}$, $X_1^j = x_{(j-1)q+1} \cdots x_{jq}$.

**Case I :** Let us assume that $r \in [0, (q-1)\,p]$. Consider the sequence $T_1 = x_{r+1} \ldots x_m$ over $C_p^d$. Then, $|T_1| = m - r \geq p\,(q+d-1) - (d-1) - p\,(q-1) = p\,d - d + 1 = D(C_p^d)$ [from 2,[O]]. Therefore, $T_1$ has a zero-sum subsequence over $C_p^d$. This implies that $T_1^S$ has a zero-sum subsequence over $G$. So, $0_G \in [S]$.

**Case II :** For $r \in [(q-1)\,p+1,\, p\,q-q]$, note that $|T| = m - r + \sum_{i=1}^{q-1} t_i \geq p\,d - d + 1 = D(C_p^d)$. Therefore, $T$ has a zero-sum subsequence $T_1$ such that $T_1^S$ has each of its $y$-coordinates $0_{C_q}$. If $T_1 = \prod_{i=1}^{(q-1)} \prod_{j=1}^{m_i} \sigma(X_i^{z_{ij}})\, x_{u_1} x_{u_2} \cdots x_{u_v}$, where $m_i \leq t_i$ and $u_i \in [r+1, m]$, then $U^S$, where $U = \prod_{i=1}^{(q-1)} \prod_{j=1}^{m_i} X_i^{z_{ij}} x_{u_1} x_{u_2} \cdots x_{u_v}$, is the required subsequence of $S$.

**Case III :** Let us assume that $r = p\,q - (q-1) + j$ where $j \in [0, q-1]$. We have, $p\,q - (q-1) + j = q \sum_{i=1}^{q-1} t_i + \sum_{i=1}^{q-1} \ell_i$, which implies $\sum_{i=1}^{q-1} \ell_i = q\,\alpha + j + 1$ for some $\alpha \in [0, q-3]$ as $\sum_{i=1}^{q-1} \ell_i \in [0, (q-1)^2]$. Note that $|T| = m - r + \sum_{i=1}^{q-1} t_i = (p\,d - d - 1) + (q - j) - \alpha$, since $\sum_{i=1}^{q-1} t_i = \frac{r - \sum_{i=1}^{q-1} \ell_i}{q} = p - 1 - \alpha$. Now consider the sequence $S_1$ of $C_q$ whose elements do not belong to the $y$-coordinate of $T^S$. Then, $|S_1| = \sum_{i=1}^{q-1} \ell_i = \alpha\,q + (j+1)$.

**Subcase III.(i)** For $j = q-1$, we have $|T| = (p\,d - d - 1) + (q - j) - \alpha = p\,d - d + 1 - (\alpha+1)$. Also, $|S_1| = (\alpha+1)\,q = (\alpha+1)\,D(C_q)$. Therefore, $S_1$ has $(\alpha+1)$ disjoint zero-sum subsequences, say, $T_k = y_{u_{k1}} \cdots y_{u_{kv_k}}$ for all $k \in [1, \alpha+1]$. Let $U_k = x_{u_{k1}} \cdots x_{u_{kv_k}}$. Assume that $T' = T\,\sigma(U_1)\,\sigma(U_2)\ldots\sigma(U_{\alpha+1})$. Then $|T'| = p\,d - d + 1 = D(C_p^d)$. Hence $0 \in [T']$. The proof follows as the corresponding $y$-coordinates of elements of $T'$ are zeros.

**Subcase III.(ii)** For $j \in [0, q-2]$ we have $|T| \geq (p\,d - d + 1) - \alpha$, and $|S_1| = \alpha\,q + j + 1 > \alpha\,q = \alpha\,D(C_q)$. So, $S_1$ has $\alpha$ disjoint zero-sum subsequences of $S_1$. Proceeding as in subcase III.(i), we conclude sub-case III.(ii). □

**Lemma 15.** *Let $G = H \times C_q$, where $q \in \mathbb{N}$ and $H$ be a finite abelian group. Let $S = (x_1, y_1) \ldots (x_n, y_n)$ be a sequence over $G$. Let $t = -\sum_{i=1}^{n} x_i$, and $\ell = -\sum_{i=1}^{n} y_i$. If $S\,(t, \ell)$ has a proper zero-sum subsequence, then $0_G \in [S]$.*

**Proof.** Let $0_G \notin [S]$, and $T$ be a proper subsequence of $S(t, \ell)$ such that $\sigma(T) = 0$. Then $T = U(t, \ell)$ for some $U \mid S$ and $|U| < |S|$. Since, $\sigma(S(t, \ell)) = (\sum_{i=1}^{n} x_i + t, \sum_{i=1}^{n} y_i + \ell) = (0_H, 0_{C_q}) = \sigma(T)$. Therefore, $\sigma(S) = \sigma(U)$ which implies $\sigma(SU^{-1}) = 0_G$, a contradiction. □

## 4 Proof of the theorem 2

**Proof of Theorem 2.** We have $(p, q) = 1$. Therefore, $G = C_p^d \times C_q \cong C_p^{d-1} \times C_{pq}$. Clearly, $D(G) \geq D^*(G) = p(q + d - 1) - d + 1$. Again, consider a sequence $S = (x_1, y_1) \cdots (x_m, y_m)$ of length $m = p(q + d - 1) - d + 1$ over $G$ such that $x_1 \cdots x_m$ is a sequence over $C_p^d$ and $y_1 \cdots y_m$ is a sequence over group $C_q$. Without loss of generality, assume that there are some $r_i \in \mathbb{N} \cup \{0\}$ for $i \in [1, q-1]$,

$$
\begin{aligned}
y_1 = \ldots = y_{r_1} &= 1, \\
y_{r_1+1} = \ldots = y_{r_1+r_2} &= 2, \\
&\vdots \\
y_{\sum_{i=1}^{q-2} r_i + 1} = \ldots = y_{\sum_{i=1}^{q-1} r_i} &= q-1, \\
y_{r+1} = \ldots = y_m &= 0 \text{ where } r = \sum_{i=1}^{q-1} r_i.
\end{aligned}
$$

For $0 \leq r \leq pq$, by Proposition 14 we can conclude that $0_G \in [S]$. On the other hand, let us assume that $pq + 1 \leq r \leq p(q + d - 1) - d + 1$. If $\sum_{i=1}^{q-1} i\, r_i \equiv 0 \pmod{q}$ then we are done by Conjecture 1. Now, assume $\sum_{i=1}^{q-1} i\, r_i \equiv j \pmod{q}$ for some $j \in [1, q-1]$ and $t = -\sum_{i=1}^{m} x_i$. Then, $(t, q - j) \in C_p^d \times C_q$. Now we have, $\sigma(S(t, q - j)) = (0_{C_p^d}, 0_{C_q})$. Consider a proper subsequence $T = S(x_m, 0)^{-1}(t, q - j)$ of $S(t, q - j)$. Clearly, $T$ is a sequence of $G$ such that $|T| = m$ and sum of the $y$-coordinates of elements of $T$ is $0_{C_q}$. Using Conjecture 1, we conclude that there exists zero-sum subsequence $T'$ of $T$. Therefore, $T'$ is proper zero-sum subsequence of $S(t, q - j)$. Now by Lemma 15, $S$ has a zero-sum subsequence. This completes the proof. □

## 5 Proof of theorem 7 and theorems related to $D_r(G)$

**Proof of Theorem 3.** For the given group $G$, if we consider the sequence $S = (1, 0, \ldots, 0)^{p^{e_1}-1}(0, 1, 0, 0, \ldots, 0)^{p^{e_2}-1} \ldots (0, 0, \ldots, 0, 1)^{rp^{e_d}-1}$, then we can have $D_r(G) \geq r\,p^{e_d} + \sum_{i=1}^{d-1} p^{e_i} - d + 1$. Also for this group $G$, $\eta(G) \leq D(G) + \exp(G)$ by Lemma 10. By Lemma 11 we have $\eta_r(G) \leq D(G) + r\exp(G)$. We conclude the proof by Lemma 12. □

**Proof of Theorem 4.**

Consider the group $G = C_{p^{e_1}} \times C_{p^{e_2}} \times \ldots \times C_{mp^{e_d}}$ such that $e_i \leq e_{i+1}$ and $p^{e_d} \geq 1 + \sum_{i=1}^{d-1} (p^{e_i} - 1)$. We need to show that

$$
D_r(G) = r\,m\,p^{e_d} + \sum_{i=1}^{d-1} p^{e_i} - d + 1.
$$

Consider the following sequence

$$S = (1, 0, \ldots, 0)^{p^{e_1}-1} (0, 1, 0, 0, \ldots, 0)^{p^{e_2}-1} \ldots (0, 0, \ldots, 1, 0)^{p^{e_{d-1}}-1} (0, 0, \ldots, 0, 1)^{rmp^{e_d}-1}.$$

Clearly the sequence $S$ does not have $r$-disjoint zero-sum subsequences. Therefore,

$$D_r(G) \geq r\,m\,p^{e_d} + \sum_{i=1}^{d-1} p^{e_i} - d + 1.$$

Let us consider $H = C_m$. Then $H$ can be seen as a normal subgroup of $G$. Then, $G/H = C_{p^{e_1}} \times C_{p^{e_2}} \times \ldots \times C_{p^{e_d}}$. By Lemma 8,

$$
\begin{aligned}
D_r(G) &\leq D_{D_r(H)}(G/H) \\
&= D_{D_r(C_m)}(C_{p^{e_1}} \times C_{p^{e_2}} \times \ldots \times C_{p^{e_d}}) \\
&= D_{rm}(C_{p^{e_1}} \times C_{p^{e_2}} \times \ldots \times C_{p^{e_d}}) \\
&= r\,m\,p^{e_d} + \sum_{i=1}^{d-1} p^{e_i} - d + 1. \quad \text{[by Theorem 3]}
\end{aligned}
$$

This completes the proof. $\square$

**Proof of Theorem 5.**

1. Consider the group $G = C_{p^{e_1}} \times C_{p^{e_2}} \times \ldots \times C_{mp^{e_{d-1}}} \times C_{np^{e_d}}$ such that $e_i \leq e_{i+1}$ and $p^{e_d} \geq 1 + \sum_{i=1}^{d-1} (p^{e_i} - 1)$.

   We need to show that if $m \,|\, n$ then,

$$r\,n\,p^{e_d} + (m-1)\,p^{e_{d-1}} + \sum_{i=1}^{d-1} p^{e_i} - d + 1 \leq D_r(G) \leq r\,n\,p^{e_d} + (m-1)\,p^{e_d} + \sum_{i=1}^{d-1} p^{e_i} - d + 1.$$

   Consider the following sequence

$$S = (1, 0, \ldots, 0)^{p^{e_1}-1} \ldots (0, \ldots, 1, 0, 0)^{p^{e_{d-2}}-1} (0, 0, \ldots, 1, 0)^{mp^{e_{d-1}}-1} (0, 0, \ldots, 0, 1)^{rnp^{e_d}-1}.$$

   Clearly the sequence $S$ does not have $r$-disjoint zero-sum subsequences. Therefore,

$$D_r(G) \geq r\,n\,p^{e_d} + (m-1)\,p^{e_{d-1}} + \sum_{i=1}^{d-1} p^{e_i} - d + 1.$$

   Let us consider $H = C_m \times C_n$. Therefore, $D_r(H) = r\,n + m - 1$ . Then, $G/H = C_{p^{e_1}} \times C_{p^{e_2}} \times \ldots \times C_{p^{e_d}}$. From Lemma 8, we have

$$
\begin{aligned}
D_r(G) &\leq D_{D_r(H)}(G/H) \\
&= D_{rn+m-1}(C_{p^{e_1}} \times C_{p^{e_2}} \times \ldots \times C_{p^{e_d}}) \\
&= (r\,n + m - 1)\,p^{e_d} + \sum_{i=1}^{d-1} p^{e_i} - d + 1 \quad \text{[by Theorem 3]}.
\end{aligned}
$$

   This completes the proof.

2. Consider the group $G = C_{p^{e_1}} \times C_{p^{e_2}} \times \ldots \times C_{mp^{e_{d-1}}} \times C_{np^{e_d}}$ such that $e_i \leq e_{i+1}$ and $p^{e_d} \geq 1 + \sum_{i=1}^{d-1} (p^{e_i} - 1)$.

We need to show that if $n \mid m$ then,

$$\max \left\{ r\, n\, p^{e_d} + (m-1)\, p^{e_d-1} + \sum_{i=1}^{d-1} p^{e_i} - d + 1, \, n\, p^{e_d} + (r\, m - 1)\, p^{e_d-1} + \sum_{i=1}^{d-1} p^{e_i} - d + 1 \right\}$$

$$\leq \quad D_r(G)$$

$$\leq \quad r\, m\, p^{e_d} + (n-1)\, p^{e_d} + \sum_{i=1}^{d-1} p^{e_i} - d + 1.$$

Consider the following two sequences

$$S_1 = (1,0,\ldots,0)^{p^{e_1}-1} \ldots (0,\ldots,1,0,0)^{p^{e_d-2}-1} (0,0,\ldots,1,0)^{m p^{e_d-1}-1} (0,0,\ldots,0,1)^{r n p^{e_d}-1}$$

and

$$S_2 = (1,0,\ldots,0)^{p^{e_1}-1} \ldots (0,\ldots,1,0,0)^{p^{e_d-2}-1} (0,0,\ldots,1,0)^{r m p^{e_d-1}-1} (0,0,\ldots,0,1)^{n p^{e_d}-1}$$

Clearly both the sequences $S_1$ and $S_2$ do not have $r$-disjoint zero-sum subsequences. Therefore,

$$D_r(G) \geq \max \left\{ r\, n\, p^{e_d} + (m-1)\, p^{e_d-1} + \sum_{i=1}^{d-1} p^{e_i} - d + 1, \, n\, p^{e_d} + (r\, m - 1)\, p^{e_d-1} + \sum_{i=1}^{d-1} p^{e_i} - d + 1 \right\}.$$

Similar to the first part, if we consider $H = C_m \times C_n$ we have, $D_r(H) = r\, m + n - 1$. Then, $G/H = C_{p^{e_1}} \times C_{p^{e_2}} \times \ldots \times C_{p^{e_d}}$. Using Lemma 8 and Theorem 3, we conclude that

$$D_r(G) \leq (r\, m + n - 1)\, p^{e_d} + \sum_{i=1}^{d-1} p^{e_i} - d + 1.$$

This completes the proof.

**Proof of Theorem 6.** Consider the group $G = C_{\prod_{j=1}^{\ell} p_j^{e_1^{(j)}}} \times \ldots \times C_{\prod_{j=1}^{\ell} p_j^{e_d^{(j)}}}$ such that $e_i^{(j)} \leq e_{i+1}^{(j)}$ and $p_j^{e_d^{(j)}} \geq 1 + \sum_{i=1}^{d-1} \left( p_j^{e_i^{(j)}} - 1 \right)$. Clearly, $D_r(G) \geq D^*(G) + (r-1)\exp(G)$. First we will prove the upper bound for $\ell = 2$ and then we proceed to prove the upper bound using inductive argument on $\ell$. For $\ell = 2$, assume that $H_1 = C_{p_1^{e_1^{(1)}}} \times \ldots \times C_{p_1^{e_d^{(1)}}}$, which implies $G/H_1 \cong C_{p_2^{e_1^{(2)}}} \times \ldots \times C_{p_2^{e_d^{(2)}}}$. Now using the Lemma 8 and the Theorem 3, we conclude that $D_r(G) \leq r\, p_1^{e_1^{(1)}} p_2^{e_d^{(2)}} + p_2^{e_2} \left( \sum_{i=1}^{d} p_1^{e_1^{(1)}} - p_1^{e_d^{(1)}} - d + 1 \right) + \left( \sum_{i=1}^{d} p_2^{e_i^{(2)}} - p_2^{e_d^{(2)}} - d + 1 \right)$. Let $\ell \geq 2$, and assume that the upper bound is true for $\ell - 1$. That is, if $H_2 = C_{\prod_{j=1}^{\ell-1} p_j^{e_1^{(j)}}} \times \ldots \times C_{\prod_{j=1}^{\ell-1} p_j^{e_d^{(j)}}}$ then $G/H_2 \cong C_{p_\ell^{e_1^{(\ell)}}} \times \ldots \times C_{p_\ell^{e_d^{(\ell)}}}$ and

$$D_r(H_2) \leq r \prod_{j=1}^{\ell-1} p_j^{e_d^{(j)}} + \sum_{m=1}^{\ell-2} \left( \left( \prod_{j=m+1}^{\ell-1} p_j^{e_d^{(j)}} \right) \varphi(p_m) \right) + \varphi(p_\ell) = k \text{ (say)},$$

where $\varphi(p_j) = \sum_{i=1}^{d-1} p_j^{e_i^{(j)}} - d + 1$ for $j \in [1, \ell - 1]$. Now by Lemma 8, we have

$$
\begin{aligned}
D_r(G) &\leq D_{D_r(H_2)} (G/H_2) \\
&\leq D_k (G/H_2) \ \left[\text{as } D_s(H_2) \leq D_{s+1}(H_2) \text{ for } s \in \mathbb{N}\right] \\
&= k\, p_d^{e_d^{(\ell)}} + \sum_{i=1}^{d-1} p_\ell^{e_i^{(\ell)}} - d + 1.
\end{aligned}
$$

This completes the proof.□

**Proof of Theorem 7.** For the group $G = C_p^{d-1} \times C_{pq}$ with $(p, q) \neq 1$, we have $q = p\, m$ for some $m \in \mathbb{N}$. Now by Lemma 8, $D^*(G) \leq D(G) \leq D_{D(C_m)} (C_p^{d-1} \times C_{p^2})$. Now as for sufficiently large prime $p$, $p^2 \geq 1 + \sum_{i=1}^{d-1} (p-1)$ and $D(C_m) = m$, then we conclude the theorem by Observation 1.□

# 6  Concluding remarks

For any finite abelian group $G$, one can observe from the theorem 6 that $D_r(G)$ must lie between the bounds that depend on $r$, primes $p_j$ and natural numbers $e_i^{(j)}$. Let us define error$:= \dfrac{\text{upper bound} - \text{lower bound}}{\text{lower bound}}$. The 'error' becomes negligible i.e., $\frac{\text{upper bound}}{\text{lower bound}} \to 1$ if either $p_j$'s are large, or $e_i^{(j)}$'s are higher natural numbers, or $r$ increases. This can be shown in general for any $\ell$ (number of primes) but then the algebra will be more complicated. Here, we will prove it for $\ell = 2$, the proof for any $\ell$ will follow similarly. For $\ell = 2$, we have

$$
\begin{aligned}
\text{error} &= \frac{\sum_{i=1}^{d-1} \left( p_1^{e_i^{(1)}} \left( p_2^{e_d^{(2)}} - p_2^{e_i^{(2)}} \right) + p_2^{e_i^{(2)}} \right) - (d-1)\, p_2^{e_d^{(2)}}}{r\, p_1^{e_d^{(1)}}\, p_2^{e_d^{(2)}} + \sum_{i=1}^{d-1} p_1^{e_i^{(1)}}\, p_2^{e_i^{(2)}} - 1} \\
&\leq \max_{i \in [1, d-1]} \frac{p_1^{e_i^{(1)}} \left( p_2^{e_d^{(2)}} - p_2^{e_i^{(2)}} \right) + \left( p_2^{e_i^{(2)}} - p_2^{e_d^{(2)}} \right)}{\dfrac{r}{d-1} \left( p_1^{e_d^{(1)}}\, p_2^{e_d^{(2)}} \right) + \left( p_1^{e_i^{(1)}}\, p_2^{e_i^{(2)}} - 1 \right)}.
\end{aligned}
$$

Clearly, the sufficient condition for 'error' to be negligible is

$$
\max_{i \in [1, d-1]} \frac{p_1^{e_i^{(1)}} \left( p_2^{e_d^{(2)}} - p_2^{e_i^{(2)}} \right) + \left( p_2^{e_i^{(2)}} - p_2^{e_d^{(2)}} \right)}{\dfrac{r}{d-1} \left( p_1^{e_d^{(1)}}\, p_2^{e_d^{(2)}} \right) + \left( p_1^{e_i^{(1)}}\, p_2^{e_i^{(2)}} - 1 \right)} \ll 1
$$

$$
\Leftrightarrow 1 - 2\, p_2^{e_i^{(2)} - e_d^{(2)}} + \frac{p_2^{e_d^{(2)}} - p_2^{e_i^{(2)}}}{p_1^{e_i^{(1)}}\, p_2^{e_d^{(2)}}} \ll \frac{r}{d-1}\, p_1^{e_d^{(1)} - e_i^{(1)}}
$$

i.e., either

1. $p_1, p_2$ are large primes, or

2. $e_i^{(1)}$'s and $e_i^{(2)}$'s are higher natural numbers, or

3. $r$ increases.

∗ Let $p$ be a prime such that $D_r(C_p^d) = (r+d-1)\,p - (d-1)$ holds for any $r$, and for some $d \in \mathbb{N}$. If $G = C_p^{d-1} \times C_{pq}$, $H = C_q$, and $r = 1$, using Lemma 8 we have $D(G) \le D_{D_r(H)}(G/H) \le D_q(C_p^{d-1}) = (q+d-1)\,p - (d-1) = D^*(G)$. This concludes that $D(G) = D^*(G)$. ⟨fi⟩